\documentclass[12pt]{article}
\usepackage{mathrsfs}
\usepackage{epic,eepic,epsf,epsfig}
\usepackage{amsfonts,srcltx,mathrsfs}
\usepackage{amsmath}
\usepackage{amssymb}
\usepackage{amsbsy}
\usepackage{graphicx}
\usepackage{amsfonts}
\usepackage{color}

\newtheorem{lem}{Lemma}[section]
\newtheorem{theorem}[lem]{Theorem}

\newtheorem{case}{Case}
\newtheorem{prop}[lem]{Proposition}

\def\a{\alpha} \def\b{\beta}  \def\d{\delta} 
 \def\s{\sigma}

\def\di{\bigm|} \def\lg{\langle} \def\rg{\rangle}
\def\nd{\mathrel{\bigm|\kern-.7em/}}

\def\f{\noindent}

\def\AGL{\hbox{\rm AGL}}

\def\Aut{\hbox{\rm Aut}}

\def\Aut{\hbox{\rm Aut}}
\def\N{\hbox{\rm N}}\def\O{\hbox{\rm O}}

\def\AGL{\hbox{\rm AGL}}

\def\ASL{\hbox{\rm ASL}}\def\AGL{\hbox{\rm AGL}}
\def\Cay{\hbox{\rm Cay}}

\def\O{\hbox{\rm O}}

\def\demo{\f {\bf Proof.}\hskip10pt}

\newcommand{\qed}{\mbox{\raisebox{0.7ex}{\fbox{}}} \vspace{4truemm}}
\def\mz{{\mathbb Z}}

\def\bin#1#2{{#1\choose#2}}

\begin{document}
\title{Half-arc-transitive graphs of prime-cube order of small valencies}

\author{ \\ Yi Wang, Yan-Quan Feng*\\
{\small\em Department of Mathematics, Beijing Jiaotong University, Beijing,
100044, P.R. China}\\}

\date{}
\maketitle

\footnotetext{*Corresponding author. E-mails: yiwang$@$bjtu.edu.cn, yqfeng$@$bjtu.edu.cn}

\begin{abstract}

A graph is called {\em half-arc-transitive} if its full automorphism group acts transitively on vertices and edges, but not on arcs. It is well known that for any prime $p$ there is no tetravalent half-arc-transitive graph of order $p$ or $p^2$. Xu~[Half-transitive graphs of prime-cube order, J. Algebraic Combin. 1 (1992) 275-282] classified half-arc-transitive graphs of order $p^3$ and valency $4$. In this paper we classify half-arc-transitive graphs of order $p^3$ and valency $6$ or $8$. In particular, the first known infinite family of half-arc-transitive Cayley graphs on non-metacyclic $p$-groups is constructed.

\bigskip
\f {\bf Keywords:} Cayley graph, half-arc-transitive graph, automorphism group.
\end{abstract}

\section{Introduction\label{s1}  }

A (di)graph $\Gamma$ consists of a pair of sets $(V(\Gamma),E(\Gamma))$, where $V(\Gamma)$ is its {\em vertex set}, and $E(\Gamma)$ is its {\em edge set}. For a graph,  $E(\Gamma)$ is also called {\em undirected edge set} that is a subset of the set $\{\{u,v\}\ |\ u,v\in V(\Gamma)\}$, and for a digraph, $E(\Gamma)$  is also called {\em directed edge set} that is a subset of the set $\{(u,v)\ |\ u,v\in V(\Gamma)\}$. For an edge $\{u,v\}$ of a graph $\Gamma$, we call $(u,v)$ an {\em arc} of $\Gamma$. An automorphism of a (di)graph $\Gamma$ is a permutation on $V(\Gamma)$ preserving the adjacency of $\Gamma$, and all automorphisms of $\Gamma$ form a group under the composition of permutations, called the {\em full automorphism group} of $\Gamma$ and denoted by $\Aut(\Gamma)$. A (di)graph $\Gamma$ is {\em vertex-transitive} or {\em edge-transitive} if $\Aut(\Gamma)$ acts transitively on $V(\Gamma)$ or $E(\Gamma)$, respectively. A graph $\Gamma$ is {\em arc-transitive} or {\em symmetric} if $\Aut(\Gamma)$ is transitive on the arc set of $\Gamma$, and {\em half-arc-transitive} provided that it is vertex-transitive, edge-transitive, but not arc-transitive. Throughout this paper, all (di)graphs $\Gamma$ are finite and simple, that is, $|V(\Gamma)|$ is finite and there is no edge between the same vertex. Let $G$ be a permutation group on a set $\Omega$ and $\alpha \in \Omega$. Denote by $G_{\alpha}$ the stabilizer of $\alpha$ in $G$, that is, the subgroup of $G$ fixing the point $\alpha$. We say that $G$ is semiregular on $\Omega$ if $G_{\alpha}=1$ for every $\alpha \in \Omega$ and regular if $G$ is transitive and semiregular.

Let $G$ be a finite group and $S$ a subset of $G$ such that $1\notin S$. The {\em Cayley  digraph} $\Gamma=\Cay(G,S)$ on $G$ with respect to $S$ is defined as the digraph with vertex set $V(\Gamma)=G$ and directed edge set $\{(g,sg)\ |\ g\in G, s\in S\}$. A Cayley digraph $\Cay(G,S)$ is connected if and only if $G=\langle S \rangle$, and if $S$ is symmetric, that is, $S^{-1}=\{ s^{-1}\ |\ s\in S\}=S$, then $\Cay(G,S)$ can be viewed as a graph by identifying the two opposite directed edges $(g,sg)$ and $(sg,g)$ as an undirected edge $\{g,sg\}$. Thus, Cayley graph can viewed as a special case of Cayley digraph. It is easy to see that $\Aut(\Cay(G,S))$ contains the right regular representation $\hat{G}=\{\hat{g}\ |\ g \in G\}$ of $G$, where $\hat{g}$ is the map on $G$ defined by $x\mapsto xg$, $x\in G$, and $\hat{G}$ is regular on the vertex set $V(\Gamma)$. This implies that a Cayley digraph is vertex-transitive. Also, it is easy to check that  $\Aut(G,S)=\{ \alpha\in \Aut(G)\ |\ S^{\alpha}=S\}$ is a subgroup of $\Aut(\Cay(G,S))_{1}$, the stabilizer of the vertex $1$ in $\Aut(\Cay(G,S))$. A Cayley digraph $\Cay(G,S)$ is said to be {\em normal} if $\hat{G}$ is normal in $\Aut(\Cay(G,S))$.

In 1947, Tutte~\cite{Tutte} initiated an investigation of half-arc-transitive graphs by showing that a vertex- and edge-transitive graph with odd valency must be arc-transitive. A few years later, in order to answer Tutte's question of the existence of half-arc-transitive graphs of even valency, Bouwer~\cite{Bouwer} gave a construction of $2k$-valent half-arc-transitive graph for every $k\geq 2$. Following these two classical articles, half-arc-transitive graphs have been extensively studied from different perspectives over decades by many authors. See, for example,~\cite{Alspach1, Alspach, Dushaofei, Holt, Li2, RJWang, Xu1}.
One of the standard problems in the study of half-arc-transitive graphs is to classify such graphs of certain orders. Let $p$ be a prime. It is well known that there are no half-arc-transitive graphs of order $p$ or $p^2$~\cite{Chao}, and by Cheng and Oxley~\cite{Cheng and Oxley}, there are no half-arc-transitive graphs of order $2p$. Alspach and Xu~\cite{Alspach} classified  half-arc-transitive graphs of order $3p$ and Wang~\cite{RJWang} classified half-arc-transitive graphs of order a product of two distinct primes. Classification of half-arc-transitive graphs of order $4p$ had been considered for more than $10$ years by many authors, and recently was solved by Kutnar et al.~\cite{Kutnar}. Despite all of these efforts, however, further classifications of half-arc-transitive graphs with general valencies seem to be very difficult.

In view of the fact that $4$ is the smallest admissible valency for a half-arc-transitive graph, special attention has rightly been given to the study of tetravalent half-arc-transitive graphs. In particular, constructing and classifying tetravalent half-arc-transitive graphs is currently an active topic in algebraic graph theory (for example, see~\cite{Feng-valency 4, Feng, Dragan-valency4, Primoz-valency4}). Maru\v{s}i\v{c}~\cite{Marusic group action} classified tightly attached tetravalent half-arc-transitive graphs with odd radius, and Maru\v{s}i\v{c} and Praeger~\cite{Marusic and Praeger} classified tightly attached tetravalent half-arc-transitive graphs with even radius.
For quite some time, all known examples of tetravalent half-arc-transitive graphs had vertex-stabilizers that are either abelian or dihedral: For instance, Maru\v{s}i\v{c}~\cite{Marusic large stabilizer} constructed an infinite family of tetravalent half-arc-transitive graphs having vertex stabilizers isomorphic to $\mathbb{Z}_{2}^{m}$ for each positive integer $m\geq 1$, and Conder and Maru\v{s}i\v{c}~\cite{Conder nonabelian stabilizer} constructed a tetravalent half-arc-transitive graph with vertex-stabilizer isomorphic to $D_4$ of order $8$.
Until recently, a tetravalent half-arc-transitive graph with vertex-stabilizers that are neither abelian nor dihedral was constructed by Conder~et~al.~\cite{Conder recent discoverry}.

Xu~\cite{Xu1} classified tetravalent half-arc-transitive graphs of order $p^3$ for each prime $p$, and later, it was extended to the case of $p^4$ by Feng et al.~\cite{Feng}. In this paper, we classify half-arc-transitive graphs of order $p^3$ and valency $6$ or $8$. In these new constructions, there is an infinite family of half-arc-transitive Cayley graphs on non-metacyclic $p$-groups, and to our best knowledge, this is the first known construction of such graphs.

From elementary group theory we know that up to isomorphism there are only five groups of order $p^3$, that is, three abelian groups $\mathbb{Z}_{p^3}$, $\mathbb{Z}_{p^2}\times \mathbb{Z}_p$ and $\mathbb{Z}_p\times \mathbb{Z}_p\times \mathbb{Z}_p$, where $\mathbb{Z}_p$ denotes the cyclic group of order $n$, and two non-abelian groups $G_1(p)$ and $G_2(p)$ defined as $$G_1(p)=\langle a,b\ |\ a^{p^2}=1, b^p=1, b^{-1}ab=a^{1+p}\rangle$$ and $$G_2(p)=\langle a,b,c\ |\ a^p=b^p=c^p=1, [a,b]=c, [a,c]=[b,c]=1 \rangle.$$

Denote by $\mathbb{Z}_{n}^{*}$ the multiplicative group of $\mathbb{Z}_n$ consisting of numbers coprime to $n$. Let $e$ be an element of order $j<p$ in $\mz_{p^2}^*$. Since $\mz_{p^2}^*\cong \mz_{p(p-1)}$, we have $j\mid (p-1)$. For each $1\leq k\leq p-1$, let
$T^{j,k}=\{b^ka,b^ka^e,\ldots,b^ka^{e^{j-1}},(b^ka)^{-1},(b^ka^e)^{-1},
\ldots,(b^ka^{e^{j-1}})^{-1}\}$ and define
$$\Gamma^{j,k}=\Cay(G_1(p),T^{j,k}).$$

Let $\lambda$ be an element of order $4$ in $\mz_p^*$. Then $4\mid (p-1)$. For each $0\leq k\leq p-1$ with $k\neq 2^{-1}(1+\lambda)$, let $S_{4,k}=\{a, b, a^\lambda b^{\lambda-1}c^k, a^{-\lambda-1}b^{-\lambda}c^{1-k}, a^{-1}, b^{-1}, (a^\lambda b^{\lambda-1}c^k)^{-1}, (a^{-\lambda-1}b^{-\lambda}c^{1-k})^{-1}\}$ and define  $$\Gamma_{4,k}=\Cay(G_2(p),S_{4,k}).$$
The following is the main result of the paper.

\begin{theorem}\label{maintheorem}
Let $\Gamma$ be a graph of order $p^3$ for a prime $p$. Then we have:
 \begin{enumerate}
   \item[\rm (1)] If $\Gamma$ has valency $6$ then $\Gamma$ is half-arc-transitive if and only if $3\mid (p-1)$ and $\Gamma\cong \Gamma^{3,k}$. There are exactly $(p-1)/2$ nonisomorphic half-arc-transitive graphs in $\Gamma^{3,k}$;
   \item[\rm (2)] If $\Gamma$ has valency $8$ then $\Gamma$ is half-arc-transitive if and only if $4\mid (p-1)$ and $\Gamma\cong \Gamma^{4,k}$ or $\Gamma_{4,k}$. There are exactly $(p-1)/2$ nonisomorphic half-arc-transitive graphs in $\Gamma^{4,k}$ and $\Gamma_{4,k}$, respectively.
 \end{enumerate}
\end{theorem}

\section{Preliminaries}\label{section=preliminaries}

We start by stating some group-theoretical results. For a group $G$ and $x, y\in G$, denote by $[x,y]$ the commutator $x^{-1}y^{-1}xy$. The following proposition is a basic property of commutators and its proof is straightforward (also see \cite[5.1.5]{Robinson}):

\begin{prop}{\rm(\cite[Kapitel III, Hilfss$\ddot{\rm a}$tze 1.2 and 1.3]{Huppert})}\label{commmutator}
Let $G$ be a group. Then, for any $x,y,z\in G$, we have $[x,y]=[y,x]^{-1}$, $[xy,z]=[x,z]^y[y,z]$ and $[x,yz]=[x,z][x,y]^z$. Furthermore, if $[x,y]$ commutes with $x$ and $y$ then for any integers $i$ and $j$, $[x^i,y^j]=[x,y]^{ij}$, and for any positive integer $n$, $(xy)^n=x^ny^n[y,x]^{\bin{n}{2}}$.
\end{prop}

In Proposition~\ref{commmutator}, it is easy to show that  the equation $(xy)^n=x^ny^n[y,x]^{\bin{n}{2}}$ holds for any integer $n$ if we set $\bin{n}{2}=\frac{n(n-1)}{2}$. The following proposition shows that there is no half-arc-transitive graphs on a group of order $p^3$ and valency $2p$.

\begin{prop}{\rm(\cite[Theorem 1.2]{Li1})} \label{prop=valency 2p symmetric}
Let $p$ be a prime, and let $\Gamma$ be a vertex-transitive and edge-transitive graph of order $p^3$ and valency $2p$. Then $\Gamma$ is symmetric.
\end{prop}

By Li and Sim~\cite[Theorem 1.1 and Lemma 2.6]{Li2}, we have the following proposition.

\begin{prop}\label{prop=metacirculant graph}
Let $\Gamma$ be a Cayley graph on $G_1(p)$ of valency $2j$ with $1<j<p$. Then $\Gamma$ is half-arc-transitive if and only if $j\big| (p-1)$ and $\Gamma\cong \Gamma^{j,k}$ for $1\leq k\leq p-1$. Furthermore,  for each $j\big| (p-1)$ there exist exactly $(p-1)/2$ nonisomorphic such graphs in $\Gamma^{j,k}$.
\end{prop}

Since a transitive permutation group of prime degree $p$ has a regular Sylow $p$-subgroup, every vertex-transitive digraph of order a prime must be a Cayley digraph. Together with the results given by $\rm Maru\check{s}i\check{c}$~\cite{Marusic}, we have the following proposition.

\begin{prop}\label{prop=vertex-transitive of $p^3$ is Cayley graph}
Any vertex-transitive digraph of order $p^k$ with $1\leq k\leq 3$ is a Cayley digraph on a group of order $p^k$.
\end{prop}

For any abelian group $H$, the map $h\mapsto h^{-1}, h\in H$ is an automorphism of $H$. By \cite[Proposition 2.10]{Feng}, we have the following proposition.

\begin{prop} \label{prop=half-arc-transitive no involution}
Let $G$ be a finite group and $\Cay(G,S)$ a connected half-arc-transitive Cayley graph. Then, $S$ does not contain an involution and for any $s\in S$, there is no $\alpha\in \Aut(G,S)$ satisfying $s^\alpha=s^{-1}$. Furthermore, every edge-transitive Cayley graph on an abelian group is also arc-transitive.
\end{prop}

\medskip
The following proposition is about isomorphisms between Cayley graphs on $p$-groups.

\begin{prop} {\rm(\cite[Theorem~1.1 (3)]{Li1})}\label{prop=CayleyIsomorphism}
Let $\Cay(G,S)$ and $\Cay(G,T)$ be two connected Cayley graphs on a $p$-group $G$ with respect to subsets $S$ and $T$, and let $|S|=|T|<2p$. Then $\Cay(G,S)$ and $\Cay(G,T)$ are isomorphic if and only if there is an automorphism $\alpha$ of $G$ such that $S^\alpha=T$.
\end{prop}

Let $X=\Cay(G,S)$ be a Cayley digraph on a finite group $G$. By Godsil~\cite[Lemma 2.2]{Godsil} (also see \cite[Proposition 1.5]{Xu}), we have $\N_{{\rm Aut}X}(\hat{G})=\hat{G}\rtimes \Aut(G,S)$.

\begin{prop}\label{prop=Cayleynormal}
A Cayley digraph $X=\Cay(G,S)$ is normal if and only if $\Aut(X)_1=\Aut(G,S)$.
\end{prop}

A finite group $G$ is called {\it $2$-genetic} if each normal subgroup of the group can be generated by two elements. The following proposition is about automorphism groups of Cayley digraphs on $2$-genetic groups.

\begin{prop}{\rm(\cite[Theorem 1.1]{Wang})}\label{3,5,7,11}
Let $G$ be a nonabelian $2$-genetic group of order $p^{n}$ for an odd prime $p$ and a positive integer $n$, and let $\Gamma=\Cay(G,S)$ be a connected Cayley digraph. Assume that $\Aut(G,S)$ is a $p'$-group and $\Gamma$ is non-normal. Then $p=3, 5, 7, 11$ and $\ASL(2,p)\leq \Aut(\Gamma)/\Phi(\O_p(A)) \leq \AGL(2,p)$. Furthermore, the kernel of $A:=\Aut(\Gamma)$ acting on the quotient digraph $\Gamma_{\Phi({\rm O}_p(A))}$ is $\Phi({\rm O}_p(A))$, and one of the following happens:
\begin{enumerate}
  \item[\rm (1)] $p=3$, $n\geq 5$, and $\Gamma_{{ \Phi({\rm O}_p(A))}}$ has out-valency at least $8$;
  \item[\rm (2)] $p=5$, $n\geq 3$ and $\Gamma_{\Phi({\rm O}_p(A))}$ has out-valency at least $24$;
  \item[\rm (3)] $p=7$, $n\geq 3$ and $\Gamma_{\Phi({\rm O}_p(A))}$ has out-valency at least $48$;
  \item[\rm (4)] $p=11$, $n\geq 3$ and $\Gamma_{\Phi({\rm O}_p(A))}$ has out-valency at least $120$.
\end{enumerate}
\end{prop}

In Proposition~\ref{3,5,7,11}, the quotient digraph $\Gamma_{\Phi({\rm O}_p(A))}$ has the orbits of $\Phi({\rm O}_p(A))$ on $V(\Gamma)$ as vertices, and for two orbits $O_1$ and $O_2$, $(O_1,O_2)$ is a directed edge in $\Gamma_{\Phi({\rm O}_p(A))}$ if and only if $(u,v)$ is a directed edge in $\Gamma$ for some $u\in O_1$ and $v\in O_2$.

\section{Proof of Theorem~\ref{maintheorem}}\label{section=proof of maintheorem}

Let $\Gamma$ be a half-arc-transitive graph and $A=\Aut(\Gamma)$. Let $(u,v)$ be an arc of $\Gamma$ and set $(u,v)^A=\{(u^a,v^a)\ |\ a\in A\}$.
Define digraphs $\Gamma_1$ and $\Gamma_2$ having vertex set $V(\Gamma)$ and directed edge sets $(u,v)^A$ and $(v,u)^A$, respectively. Since $\Gamma$ is half-arc-transitive, for every edge $\{x,y\}\in E(\Gamma)$, both $\Gamma_1$ and $\Gamma_2$ contains exactly one of the directed edges $(x,y)$ and $(y,x)$, and $\Gamma$ is connected if and only if $\Gamma_i$ is connected for each $i=1,2$. Furthermore, $A=\Aut(\Gamma_i)$  and $\Gamma_i$ is $A$-edge-transitive. In what follows we denote by $\overrightarrow{\Gamma}$ one of the digraphs
$\Gamma_1$ and $\Gamma_2$.

Let $\Gamma$ be a half-arc-transitive graph of order $p^3$ for a prime $p$.
Since there exists no half-arc-transitive graph of order less than $27$ (see~\cite{Alspach1}),  we have $p\geq 3$. By Proposition~\ref{prop=vertex-transitive of $p^3$ is Cayley graph}, $\Gamma=\Cay(G,S)$ and $\overrightarrow{\Gamma}=\Cay(G,R)$. Since a group of order $p$ or $p^2$ is abelian and there is no half-arc-transitive Cayley graph on an abelian group by Proposition~\ref{prop=half-arc-transitive no involution}, $\Gamma$ is connected, and so $G=\lg R\rg$ and $S=R\cup R^{-1}$. Furthermore, $G=G_1(p)$ or $G_2(p)$, where

\vskip 0.2cm
$\begin{array}{rl}
  G_1(p)= & \langle a,b\ |\ a^{p^2}=1, b^p=1, b^{-1}ab=a^{1+p}\rangle,  \\
  G_2(p)= & \langle a,b,c\ |\ a^p=b^p=c^p=1, [a,b]=c, [a,c]=[b,c]=1 \rangle.
\end{array}$

\vskip 0.2cm
Since $G=\langle R\rangle$ is non-abelian, $R$ contains two elements $x$ and $y$ such that $xy\not=yx$, and since $|G|=p^3$, we have $\langle x,y\rangle=G$.
For $G=G_2(p)$, $x$ and $y$  have the same relations as do $a$ and $b$, which implies that $x\mapsto a$ and $y\mapsto b$ induce an automorphism of $G_2(p)$. In this case, we may assume that $a,b\in R$. Similarly, for $G=G_1(p)$ we may assume that $a\in R$. Thus we have the following Observation:

\vskip 0.2cm

\f {\bf Observation:} {\it Let $\Gamma$ be a half-arc-transitive graph of order $p^3$ for a prime $p$. Then $\Gamma=\Cay(G,S)$ and $\overrightarrow{\Gamma}=\Cay(G,R)$ for some group $G=G_1(p)$ or $G_2(p)$ with $p\geq 3$, where $G=\lg R\rg$ and $S=R\cup R^{-1}$. Furthermore,
\begin{enumerate}
\item[\rm (1)] if $G=G_1(p)$ then $a\in R$;
\item[\rm (2)] if  $G=G_2(p)$ then $a,b\in R$.
\end{enumerate}}

Let us begin by considering normal half-arc-transitive Cayley graphs on $G_2(p)$.

\medskip
\begin{lem}\label{lemma=necessity valency 8}
Let $\Gamma=\Cay(G_2(p),S)$ be a Cayley graph of valency $8$. Then $\Gamma$ is normal and half-arc-transitive if and only if $4\big| (p-1)$ and $\Gamma\cong \Gamma_{4,k}$. Furthermore, there are exactly $(p-1)/2$ nonisomorphic half-arc-transitive graphs in $\Gamma_{4,k}$.
\end{lem}

\demo Let $\Gamma=\Cay(G_2(p),S)$ be normal and half-arc-transitive. Set $A=\Aut(\Gamma)$. By Observation, we have  $\overrightarrow{\Gamma}=\Cay(G_2(p),R)$ with $p\geq 3$, $G_2(p)=\langle R\rangle$ and $S=R\cup R^{-1}$. We may further assume $a,b,a^ib^jc^k\in R$. Since $\Gamma$ has valency $8$, we have $|S|=8$ and $|R|=4$. Since $\Gamma$ is normal, Proposition~\ref{prop=Cayleynormal} implies that $A_1=\Aut(G_2(p),S)=\Aut(G_2(p),R)$, which is transitive on $R$. Since $|R|=4$, $\Aut(G_2(p),R)\leq S_4$. It follows that  either $\Aut(G_2(p),R)$ contains two involution $\a_1$ and $\a_2$ such that $\langle \a_1,\a_2\rangle\cong\mz_2\times\mz_2$ is regular on $R$, or $\Aut(G_2(p),R)$ has an element $\a$ of order $4$ such that $\langle \a\rangle\cong\mz_4$ is regular on $R$.

Suppose that $\Aut(G_2(p),R)$ contains two involution $\a_1$ and $\a_2$ such that $\langle \a_1,\a_2\rangle\cong\mz_2\times\mz_2$ is regular on $R$.
Without loss of generality, we may assume that $a^{\a_1}=b$ and $b^{\a_1}=a$, and so $c^{\a_1}=c^{-1}$. This yields that $R=\{a,b,a^ib^jc^k,(a^ib^jc^k)^{\a_1}\}=\{a,b,a^ib^jc^k,a^jb^ic^{-ij-k}\}$.
Since $\langle \a_1,\a_2\rangle\cong\mz_2\times\mz_2$, we may assume that $a^{\a_2}=a^ib^jc^k$ and $(a^ib^jc^k)^{\a_2}=a$. Then $b^{\a_2}=a^jb^ic^{-ij-k}$ and so $c^{\a_2}=c^{i^2-j^2}$.
By Proposition~\ref{commmutator}, $a=a^{\alpha_{2}^{2}}=(a^ib^jc^k)^{\a_2}=a^{i^2+j^2}b^{2ij}c^{-ij^3+k(i^2-j^2+i-j)-ij^2-2^{-1}i^2j(i-1)-2^{-1}ij^2(j-1)}$, implying the following equations:

\parbox{8cm}{
\begin{eqnarray*}
&& i^2+j^2=1;\\
&& 2ij=0;\\
&& -ij^3+k(i^2-j^2+i-j)-ij^2-2^{-1}i^2j(i-1)-2^{-1}ij^2(j-1)=0.
\end{eqnarray*}}\hfill
\parbox{1cm}{
\begin{eqnarray}
\label{eq4}\\ \label{eq5}\\\label{eq6}
\end{eqnarray}}

\f As above, in what follows all equations are considered in $\mathbb{Z}_p$, unless otherwise stated. By Eq~(\ref{eq5}), we have $i=0$ or $j=0$.
For $i=0$, we have $j^2=1$ by Eq~(\ref{eq4}), that is, $j=\pm1$. If $j=-1$ then $S=\{a, b, a^{-1}c^{-k}, b^{-1}c^k\}\cup \{a^{-1}, b^{-1}, ac^{k}, bc^{-k}\}$, and the automorphism of $G_2(p)$ induced by $a \mapsto a^{-1}$, $b \mapsto bc^{-k}$, $c \mapsto c^{-1}$, fixes $S$ setwise, contrary to Proposition~\ref{prop=half-arc-transitive no involution}. If $j=1$ then $k=0$ by Eq~(\ref{eq6}), implying that $a^ib^jc^k=b$, a contradiction. For $j=0$, we have $i^2=1$ by Eq~(\ref{eq4}), that is, $i=\pm1$. If $i=-1$ then $S=\{a, b, a^{-1}c^k, b^{-1}c^{-k}\}\cup \{a^{-1}, b^{-1}, ac^{-k}, bc^{k}\}$, and the automorphism of $G_2(p)$ induced by $a \mapsto a^{-1}$, $b \mapsto bc^{k}$, $c \mapsto c^{-1}$, fixes $S$ setwise, contrary to Proposition~\ref{prop=half-arc-transitive no involution}. If $i=1$, then $k=0$ by Eq~(\ref{eq6}), implying that $a^ib^jc^k=a$, a contradiction.

Thus, $\Aut(G_2(p),R)$ has an element $\a$ of order $4$ and $\langle \a\rangle\cong\mz_4$ is regular on $R$.
Then $R=\{a, a^{\a}, a^{\a^2}, a^{\a^3}\}$, and since $G_2(p)=\langle R\rangle$, we have $\langle a,a^{\a}\rangle=G_2(p)$ and so $a\mapsto a$, $b\mapsto a^{\a}$ induces an automorphism of $G_2(p)$. This implies that we may assume that $a^{\a}=b$, and $\a$ is induced by $a\mapsto b$, $b\mapsto a^ib^jc^k$, $c\mapsto c^{-i}$. It follows that $R=\{a, b, a^ib^jc^k, a^{ij}b^{i+j^2}c^{-i^2j+k(j-i)-2^{-1}ij^2(j-1)}\}$.

Since  $a=a^{\alpha^4}
=a^{i(i+j^2)}b^{j(2i+j^2)}c^{i^3j+(k-i^2j)(i+j^2)-ik(j-i)+2^{-1}i^2j^2(j-1)-2^{-1}ij(i+j^2)(i+j^2-1)}$, we have the following equations:

\f \parbox{15cm}{
\begin{eqnarray*}
&&i(i+j^2)=1;\\
&&j(2i+j^2)=0;\\
&&i^3j+(k-i^2j)(i+j^2)-ik(j-i)+2^{-1}i^2j^2(j-1)-2^{-1}ij(i+j^2)(i+j^2-1)=0.
\end{eqnarray*}}\hfill
\parbox{1cm}{
\begin{eqnarray}
\label{eq7}\\ \label{eq8}\\\label{eq9}
\end{eqnarray}}

\f By Eq~(\ref{eq8}), either $j=0$ or $2i+j^2=0$. Suppose $j=0$. By Eq~(\ref{eq7}), $i^2=1$ , that is, $i=\pm1$. If $i=1$ then $k=0$ by Eq~(\ref{eq9}), and hence $a^ib^jc^k=a$, a contradiction. If $i=-1$ then $S=\{a, b, a^{-1}c^k, b^{-1}c^k\}\cup\{a^{-1}, b^{-1}, ac^{-k}, bc^{-k}\}$, and the automorphism of $G_2(p)$ induced by $a \mapsto a^{-1}$, $b \mapsto bc^{-k}$, $c \mapsto c^{-1}$, fixes $S$ setwise, contrary to Proposition~\ref{prop=half-arc-transitive no involution}.
Thus, $2i+j^2=0$. By Eq~(\ref{eq7}), $i^2=-1$ and so $ij^2=1-i^2=2$.

Thus, $S=\{a, b, a^ib^jc^{k}, a^{ij}b^{-i}c^{1+k(j-i)}\}\cup \{a^{-1}, b^{-1}, a^{-i}b^{-j}c^{-k-ij}, a^{-ij}b^ic^{-1-k(j-i)-j}\}$ because  $-i^2j+k(j-i)-2^{-1}ij^2(j-1)=1+k(j-i)$. Noting that  $i^2=-1$, $ij^2=2$ and $i+j^2=-i$, Eq~(\ref{eq9}) implies that  $2k(1+i+ij)=j-2i-ij$, and since $j^2=-2i$, we have $2kj(1+i+ij)=j(j-2i-ij)=-2(1+i+ij)$. It follows that $(kj+1)(1+i+ij)=0$, that is, either $kj+1=0$ or $1+i+ij=0$.

For $kj+1=0$, we have $k=-j^{-1}$ and so $1+k(j-i)=ij^{-1}$. It follows that
$S=\{a, b, a^ib^jc^{-j^{-1}}, a^{ij}b^{-i}c^{ij^{-1}}\}\cup \{a^{-1}, b^{-1}, a^{-i}b^{-j}c^{j^{-1}-ij}, a^{-ij}b^ic^{-ij^{-1}-j}\}$.
Note that $ij^2=2$ implies that  $-j^{-1}=j^{-1}-ij$, and $i+j^2=-i$ implies that $ij^{-1}=-ij^{-1}-j$. Then  the automorphism of $G_2(p)$ induced by $a \mapsto a^{-1}$, $b \mapsto b^{-1}$, $c \mapsto c$, fixes $S$ setwise, contrary to Proposition~\ref{prop=half-arc-transitive no involution}.

For $1+i+ij=0$, we have $0=i(1+i+ij)=i-1-j$. It follows that $j=i-1$ and
$S=\{a, b, a^ib^{i-1}c^k, a^{-i-1}b^{-i}c^{1-k}\}\cup \{a^{-1}, b^{-1}, a^{-i}b^{1-i}c^{-k+1+i}, a^{i+1}b^ic^{k-i}\}$. If $k=2^{-1}(1+i)$ then $k=-k+1+i$ and $1-k=k-i$, and then the automorphism of $G_2(p)$ induced by $a\mapsto a^{-1}$, $b\mapsto b^{-1}$, $c\mapsto c$, fixes $S$ setwise, contrary to Proposition~\ref{prop=half-arc-transitive no involution}. Hence $k\not=2^{-1}(1+i)$. Note that $i^2=-1$ implies that $4\di (p-1)$ and $i=\lambda$ is an element of order $4$ in $\mz_p^*$. Then $j=\lambda-1$ and  $k\not=2^{-1}(1+\lambda)$. By the definition of $\Gamma_{4,k}$ before Theorem~\ref{maintheorem},
$\Gamma\cong \Gamma_{4,k}$.

Now we prove that $\Gamma_{4,k}=\Cay(G_2(p),S_{4,k})$ is normal and half-arc-transitive. Recall that $S_{4,k}=\{a, a^{-1}, b, b^{-1}, a^\lambda b^{\lambda-1}c^k, a^{-\lambda}b^{1-\lambda}c^{-k+\lambda+1}, a^{-\lambda-1}b^{-\lambda}c^{1-k}, a^{\lambda+1}b^\lambda c^{k-\lambda}\}$, $\lambda$ is an element of order $4$ in $\mz_p^*$, and $k\neq 2^{-1}(1+\lambda)$. Let $A=\Aut(\Gamma_{4,k})$ and $R_{4,k}=\{a, b, a^\lambda b^{\lambda-1}c^k, a^{-\lambda-1}b^{-\lambda}c^{1-k}\}$. Then $S_{4,k}=R_{4,k}\cup R_{4,k}^{-1}$.
Let $\alpha$ be the automorphism of $G_2(p)$ induced by
$a\mapsto b$, $b\mapsto a^{\lambda}b^{{\lambda}-1}c^k$, $c\mapsto c^{-{\lambda}}$. By Proposition~\ref{commmutator},  $(a^{\lambda}b^{{\lambda}-1}c^k)^\alpha=a^{-{\lambda}-1}b^{-{\lambda}}c^{1-k}$ and $(a^{-{\lambda}-1}b^{-{\lambda}}c^{1-k})^\alpha=a$. Thus, $\alpha\in \Aut(G_2(p),S_{4,k})$ has order $4$ and permutes the elements of $R_{4,k}$ cyclically, implying that $\hat{G}_2(p)\rtimes\langle\a\rangle$ is half-arc-transitive on $\Gamma_{4,k}$. To prove the normality and the half-arc-transitivity of $\Gamma_{4,k}$, it suffices to show that $A=\hat{G}_2(p)\rtimes\langle\a\rangle$.

Write $L=\Aut(G_2(p),S_{4,k})$. Then $L$ acts on $S_{4,k}$ faithfully. Set $\Omega_1=\{a,a^{-1}\}$, $\Omega_2=\{b,b^{-1}\}$, $\Omega_3=\{a^\lambda b^{\lambda-1}c^k,a^{-\lambda}b^{1-\lambda}c^{-k+\lambda+1}\}$ and $\Omega_4=\{a^{-\lambda-1}b^{-\lambda}c^{1-k}, a^{\lambda+1}b^\lambda c^{k-\lambda}\}$.
Since $L\leq \Aut(G_2(p))$, $\{\Omega_1, \Omega_2, \Omega_3, \Omega_4\}$ is a complete imprimitive block system of $L$ on $S_{4,k}$. Let $\Omega=\{\Omega_1, \Omega_2, \Omega_3, \Omega_4\}$. Since $\a\in L$, $L$ is transitive on $\Omega$.

Let $\b\in L_a$. Then $a^\b=a$ and $\Omega_1^\b=\Omega_1$. Thus, $(\Omega_2\cup\Omega_3\cup\Omega_4)^\b=\Omega_2\cup\Omega_3\cup\Omega_4$,
and so $b^\b\in \Omega_2\cup\Omega_3\cup\Omega_4$, that is, $b^\b=b$, $b^{-1}$, $a^\lambda b^{\lambda-1}c^k$, $a^{-\lambda}b^{1-\lambda}c^{-k+\lambda+1}$, $a^{-\lambda-1}b^{-\lambda}c^{1-k}$ or $a^{\lambda+1}b^\lambda c^{k-\lambda}$.
As $\lambda$ is an element of order $4$ in $\mz_p^*$, we have $\lambda\neq 0$, $\pm1$.
If $b^\b=b^{-1}\in \Omega_2$ then $c^\b=c^{-1}$ and $\Omega_2^\b=\Omega_2$. It follows that $(\Omega_3\cup\Omega_4)^\b=\Omega_3\cup\Omega_4$, implying that $(a^\lambda b^{\lambda-1}c^k)^\b=a^\lambda b^{1-\lambda}c^{-k}\in \Omega_3\cup \Omega_4$, which is impossible.
If $b^\b= a^\lambda b^{\lambda-1}c^k\in \Omega_3$ then $c^\b=c^{\lambda-1}$ and $\Omega_2^\b=\Omega_3$. Thus, $\Omega_3^\b\subseteq \Omega_2\cup\Omega_4$, but $(a^\lambda b^{\lambda-1}c^k)^\b=a^{-1}b^{-2\lambda}c^{-\lambda+2+2k(\lambda-1)}\not\in \Omega_2\cup\Omega_4$, a contradiction.
If $b^\b=a^{-\lambda}b^{1-\lambda}c^{-k+\lambda+1}\in \Omega_3$, then $c^\b=c^{1-\lambda}$ and $\Omega_2^\b=\Omega_3$. Thus, $\Omega_3^\b\subseteq \Omega_2\cup\Omega_4$ and $(a^\lambda b^{\lambda-1}c^k)^\b=a^{2\lambda+1}b^{2\lambda}c^{-\lambda-2k(\lambda-1)}$ implies that $(a^\lambda b^{\lambda-1}c^k)^\b=b^{-1}$. It follows that $\Omega_3^\b=\Omega_2$ and so $\Omega_4^\b=\Omega_4$, which is impossible because $(a^{\lambda+1}b^{\lambda}c^{k-\lambda})^\b=$ $a^{\lambda+2}b^{\lambda+1}c^{-2k\lambda+k-\lambda-2}\not\in \Omega_4$.
If $b^\b=a^{-\lambda-1}b^{-\lambda}c^{1-k}\in \Omega_4$, then $c^\b=c^{-\lambda}$ and $\Omega_2^\b=\Omega_4$. Thus, $\Omega_3^\b\subseteq \Omega_2\cup\Omega_3$ and
$(a^\lambda b^{\lambda-1}c^{k})^\b=a^{\lambda+2}b^{\lambda+1}c^{-2k\lambda+k-\lambda-2}$ implies that $(a^\lambda b^{\lambda-1}c^{k})^\b=b^{-1}$ and $\lambda=-2$. It follows that $\Omega_3^\b=\Omega_2$ and so $\Omega_4^\b=\Omega_3$. This forces $\lambda=2$ as $(a^{\lambda+1}b^{\lambda}c^{k-\lambda})^\b=a^2bc^{-2k\lambda+\lambda-2}\in \Omega_3$, and hence $\lambda=2=-2$, a contradiction.
If $b^\b=a^{\lambda+1}b^{\lambda}c^{k-\lambda}\in \Omega_4$ then $c^\b=c^{\lambda}$ and $\Omega_2^\b=\Omega_4$. Thus, $\Omega_3^\b\subseteq \Omega_2\cup\Omega_3$ and so $(a^\lambda b^{\lambda-1}c^k)^\b=a^{\lambda-2}b^{-\lambda-1}c^{2k\lambda-k-\lambda}\in \Omega_2\cup\Omega_3$, which is impossible. The above arguments mean that $b^\b=b$, implying $\b=1$. Thus, $L_a=1$, and so $L_{a^{-1}}=1$ as $L\leq \Aut(G_2(p))$. Since $\a$ is transitive on $\Omega$, $L_x=1$ for any  $x\in S_{4,k}$.

Let $K$ be the kernel of $L$ on $\Omega$. Then $K$ fixes each $\Omega_i$ setwise, and hence $|K|=|K_a||a^K|\leq 2$. Suppose that $|K|=2$. Then the unique involution, say $\gamma$, in $K$ interchanges the two elements in each $\Omega_i$ because $L_x=1$. In particular, $\gamma$ is induced by
$a^\gamma=a^{-1}$, $b^\gamma=b^{-1}$ and $c^\gamma=c$.
It follows that $(a^\lambda b^{\lambda-1}c^k)^\gamma=a^{-\lambda}b^{1-\lambda}c^k$, and since  $(a^\lambda b^{\lambda-1}c^k)^\gamma\in \Omega_3$, we have $a^{-\lambda}b^{1-\lambda}c^k=a^{-\lambda}b^{1-\lambda}c^{-k+\lambda+1}$, forcing that $k=2^{-1}(1+\lambda)$, a contradiction. Thus, $K=1$ and $L\leq S_4$, the symmetric group of degree $4$.

Since $L_x=1$ for any $x\in S_{4,k}$, $L$ is semiregular on $S_{4,k}$, and so $|L|$ is a divisor of $8$. Since $\alpha\in L$, we have $|L|=4$ or $8$. Suppose $|L|=8$. Since $L\leq S_4$, $L$ is the dihedral group of order $8$, and so $\a^2\in Z(L)$. Note that $\a^2$ interchanges $\Omega_1$ and $\Omega_3$, and $\Omega_2$ and $\Omega_4$. Then $L_{\Omega_1}=L_{\Omega_1}^{\a^2}=L_{\Omega_3}$.
Since $L$ is transitive on $\Omega$, $|L_{\Omega_1}|=2$. Let $\d$ be the unique involution in $L_{\Omega_1}$. Then $\Omega_1^\d=\Omega_1$ and $\Omega_3^\d=\Omega_3$. Since $L_a=1$, we have
$a^\d=a^{-1}$, and since $K=1$, we have  $\Omega_2^\d=\Omega_4$. On the other hand, $\langle \a\rangle\unlhd L$ and so $R_{4,k}$ is an imprimitive block of $L$, yielding $R_{4,k}^\d=R_{4,k}^{-1}$. It follows that $b^\d\in \Omega_2^\d\cap R_{4,k}^\d=\Omega_4\cap R_{4,k}^{-1}$, that is, $b^\d=a^{\lambda+1}b^\lambda c^{k-\lambda}$.
Thus, $c^\delta=c^{-\lambda}$, and since $\d$ is an involution, $b=(a^{\lambda+1}b^\lambda c^{k-\lambda})^\d=a^{-2}b^{-1}c^{-1}$, which is impossible.
Thus, $|L|=4$ and $L=\langle \a\rangle$. Clearly, $p\nmid |L|=|\Aut(G_2(p),S_{4,k})|$. By Proposition~\ref{3,5,7,11}, $\Gamma_{4,k}$ is a normal Cayley graph, and by Proposition~\ref{prop=Cayleynormal}, $A=\hat{G}_2(p)\rtimes\langle\a\rangle$, as required.

At last, we prove the second part of the lemma by determining the number of nonisomorphic half-arc-transitive graphs in $\Gamma_{4,k}$.

There are exactly two elements of order $4$ in $\mathbb{Z}_{p}^{*}$, that is, $\lambda$ and $\lambda^{-1}=-\lambda$. Replacing $\lambda$ by $-\lambda$ in $R_{4,k}$ and $S_{4,k}$, we set $\overline{R}_{4,s}=\{a,b,a^{-\lambda}b^{-\lambda-1}c^{s}, a^{\lambda-1}b^{\lambda}c^{1-s}\}$ and $\overline{S}_{4,s}=\overline{R}_{4,s}\cup \overline{R}_{4,s}^{-1}$, where $s\in \mathbb{Z}_p$ and $s\neq 2^{-1}(1-\lambda)$. For each $k\in \mathbb{Z}_p$ and $k\neq 2^{-1}(1+\lambda)$, the automorphism $\phi$ of $G_2(p)$ induced by $a\mapsto a$, $b\mapsto a^{\lambda-1}b^\lambda c^{1-k+\lambda}$, $c\mapsto c^{\lambda}$, maps $R_{4,k}$ to $\overline{R}_{4,k-\lambda}$. Then $S_{4,k}^\phi=\overline{S}_{4,k-\lambda}$, and so $\Cay(G_2(p),S_{4,k})\cong \Cay(G_2(p),\overline{S}_{4,k-\lambda})$. Since $k\neq 2^{-1}(1+\lambda)$, we have $k-\lambda\neq 2^{-1}(1-\lambda)$. Thus, $\Gamma_{4,k}$ is independent of the choice of $\lambda$, and there are at most $p-1$ nonisomorphic half-arc-transitive graphs in $\Gamma_{4,k}$ ($k\neq 2^{-1}(1+\lambda)$). To prove there are exactly $(p-1)/2$ nonisomorphic half-arc-transitive graphs in $\Gamma_{4,k}$, it suffices to show that $\Gamma_{4,k}\cong \Gamma_{4,l}$ ($k,l\not=2^{-1}(1+\lambda))$ if and only if $l=k$ or $l=1+\lambda-k$.

Let $l=1+\lambda-k$. One may easily show that the automorphism of $G_2(p)$ induced by
$a\mapsto a^{-1}$, $b\mapsto b^{-1}$, $c\mapsto c$, maps $R_{4,k}$ to $R_{4,1+\lambda-k}^{-1}$, and so $\Gamma_{4,k}\cong \Gamma_{4,1+\lambda-k}$.

Let $\Gamma_{4,k}\cong \Gamma_{4,l}$ ($k,l\not=2^{-1}(1+\lambda))$. The stabilizers of $1$ in $\Aut(\Gamma_{4,k})$ and $\Aut(\Gamma_{4,l})$ are  $\Aut(G_2(p),S_{4,k})$ and $\Aut(G_2(p),S_{4,l})$ respectively, of which both are isomorphic to $\mz_4$.
Furthermore, $\Aut(G_2(p),S_{4,k})$ is transitive on both $R_{4,k}$ and $R_{4,k}^{-1}$, and $\Aut(G_2(p),S_{4,l})$ is transitive on both $R_{4,l}$ and $R_{4,l}^{-1}$. Since $4\di (p-1)$, we have $p\geq 5$, and by  Proposition~\ref{prop=CayleyIsomorphism}, there is an automorphism $\sigma$ of $G_2(p)$ such that $S_{4,k}^\sigma=S_{4,l}$. Thus, $\Aut(G_2(p),S_{4,k})^\sigma=\Aut(G_2(p),S_{4,l})$, and so $R_{4,k}^\sigma=R_{4,l}$ or $R_{4,l}^{-1}$. Since $\Aut(\Gamma_{4,l})$ is transitive on both $R_{4,l}$ and $R_{4,l}^{-1}$, we may assume that if $R_{4,k}^\sigma=R_{4,l}$ then $a^\sigma=a$ and if $R_{4,k}^\sigma=R_{4,l}^{-1}$ then $a^\sigma=a^{-1}$.

Assume $R_{4,k}^\sigma=R_{4,l}$ and $a^\sigma =a$. Then $b^\sigma\in R_{4,l}$ and so $b^\sigma=b$, $a^\lambda b^{\lambda-1}c^l$ or $ a^{-\lambda-1}b^{-\lambda}c^{1-l}$. If $b^\sigma=a^\lambda b^{\lambda-1}c^l$ then $c^\sigma= c^{\lambda-1}$, and hence $(a^\lambda b^{\lambda-1}c^k)^\sigma=a^{-1}b^{-2\lambda}c^{-\lambda+2+(k+l)(\lambda-1)}\in R_{4,l}$, which is impossible. If $b^\sigma= a^{-\lambda-1}b^{-\lambda}c^{1-l}$ then $c^\sigma=c^{-\lambda}$, and hence $(a^\lambda b^{\lambda-1}c^k)^\sigma=a^{\lambda+2}b^{\lambda+1}c^{-\lambda-2-l(\lambda-1)-k\lambda}\in R_{4,l}$, which is impossible. If $b^\sigma =b$ then $c^\sigma=c$, and hence $(a^\lambda b^{\lambda-1}c^k)^\sigma=a^\lambda b^{\lambda-1}c^k\in R_{4,l}$, implying that $l=k$.

Assume $R_{4,k}^\sigma=R_{4,l}^{-1}$ and  $a^\sigma=a^{-1}$. Then $b^\sigma\in R_{4,l}^{-1}$ and so $b^\sigma=b^{-1}$, $a^{-\lambda}b^{-\lambda+1}c^{-l+\lambda+1}$ or $a^{\lambda+1}b^\lambda c^{l-\lambda}$. If $b^\sigma= a^{-\lambda}b^{1-\lambda}c^{-l+\lambda+1}$ then $c^\sigma=c^{\lambda-1}$ and $(a^\lambda b^{\lambda-1}c^k)^\sigma=ab^{2\lambda}c^{-\lambda+(k-l)(\lambda-1)}\in R_{4,l}^{-1}$, which is impossible. If $b^\sigma= a^{\lambda+1}b^\lambda c^{l-\lambda}$ then $c^\sigma=c^{-\lambda}$ and hence $(a^\lambda b^{\lambda-1}c^k)^\sigma=a^{-\lambda-2}b^{-\lambda-1}c^{-\lambda-k\lambda+l(\lambda-1)}\in R_{4,l}^{-1}$, which is impossible. If $b^\sigma= b^{-1}$ then we have $c^\sigma= c$ and $(a^\lambda b^{\lambda-1}c^k)^\sigma=a^{-\lambda}b^{1-\lambda}c^k\in R_{4,l}^{-1}$, implying that $l=1+\lambda-k$. \hfill\qed

By Xu~\cite{Xu1}, there is a unique half-arc-transitive graph of order $27$ and valency $4$. However, there is no such graph of valency $6$ or $8$.

\begin{lem}\label{lem=no half 27 verties}
There is no half-arc-transitive graph of order $27$ and valency $6$ or $8$.
\end{lem}

\demo  By Proposition~\ref{prop=valency 2p symmetric}, there is no half-arc-transitive graph of order $27$ and valency $6$. Suppose that $\Gamma$ is a half-arc-transitive graph of order $27$ and valency $8$. To finish the proof, we aim to find contradictions.
By Observation, $\Gamma=\Cay(G,S)$ and $\overrightarrow{\Gamma}=\Cay(G,R)$ for $G=G_1(3)$ or $G_2(3)$, where $G=\langle R\rangle$ and  $S=R\cup R^{-1}$ with $|R|=4$. If $3\nmid |\Aut(G,S)|$ then by Proposition~\ref{3,5,7,11}, $\Gamma$ and $\overrightarrow{\Gamma}$ are normal. By Proposition~\ref{prop=Cayleynormal}, $A_1=\Aut(G,S)=\Aut(G,R)$ and hence $\Aut(G,R)$ is transitive on $R$. It follows that $4\di |\Aut(G)|$, which is impossible for $G=G_1(3)$ by \cite[Lemma 2.8]{Xu1}. Thus, $G=G_2(3)$, but this is also impossible by Lemma~\ref{lemma=necessity valency 8}. It follows that  $3\di |\Aut(G,S)|$. Let $\sigma\in \Aut(G,R)$ be of order $3$. Then $\sigma$ fixes one vertex and has an orbit of length $3$ in $R$.

\medskip
\f{\bf Case 1:} $G=G_1(3)=\langle a,b\ |\ a^{9}=1, b^3=1, b^{-1}ab=a^{1+3}\rangle$.

Since $\sigma$ has order $3$, it is induced by $a\mapsto b^sa^{1+3t}$, $b\mapsto ba^{3r}$, where $s,t,r=0,1,2$.

Assume that the fixed element of $\sigma$ in $R$ is of  order $9$. Noting that $\Aut(G)$ is transitive on elements of order $9$ in $G$, we may assume that $a\in R$ and $a^\sigma=a$. In this case, $\s$ is induced by $a\mapsto a,\ b\mapsto ba^{3r}$, where $3r=3$ or $6$.  Let $b^ia^j\in R$ with $a\not=b^ia^j$. Since $[a,b]=a^3\in Z(G)$, by Proposition~\ref{commmutator}, we have $(b^ia^j)^\s=b^ia^{j+3ri}$ and $(b^ia^j)^{\s^2}=b^ia^{j+6ri}$. It follows that $R=\{a,b^ia^j,(b^ia^j)^\s, (b^ia^j)^{\s^2}\}=\{a, b^ia^j, b^ia^{j+3ri}, b^ia^{j+6ri}\}=\{a, b^ia^j, b^ia^{j+3}, b^ia^{j+6}\}$ because $G=\langle R\rangle$ implies that $i=1$ or $2$. For $i=1$, we have $S=\{a, ba^j, ba^{j+3}, ba^{j+6}\}\cup \{a^{-1}, b^{-1}a^{2j}, b^{-1}a^{2j+3}, b^{-1}a^{2j+6}\}$, and the automorphism of $G$ induced by $a\mapsto a^{-1}$, $b\mapsto ba^{2j}$, fixes $S$ setwise, contrary to Proposition~\ref{prop=half-arc-transitive no involution}. For $i=2$, we have $S=\{a, b^2a^j, b^2a^{j+3}, b^2a^{j+6}\}\cup \{a^{-1}, ba^{2j}, ba^{2j+3}, ba^{2j+6}\}$, which is also impossible because the automorphism of $G$ induced by $a\mapsto a^{-1}$, $b\mapsto ba^j$, fixes $S$ setwise.

Assume that the fixed element of $\sigma$ in $R$ is not of order $9$. Then it has order $3$. By Proposition~\ref{commmutator}, $(b^ia^j)^3=1$ if and only if $3\di j$. Thus, we may let $b^ia^{3j}\in R$ and $(b^ia^{3j})^\sigma=b^ia^{3j}$. By Observation, $a\in R$, and so $R=\{a, a^\sigma, a^{\sigma^2},  b^ia^{3j}\}=\{a, b^sa^{1+3t}, b^{2s}a^{1+6t+3rs}, b^ia^{3j}\}$.

Let $i=0$. Then $R=\{a, b^sa^{1+3t}, b^{2s}a^{1+6t+3rs}, a^{3j}\}$, and $s\neq 0$ as $\langle R\rangle=G$. Thus, $S=\{a, b^sa^{1+3t}, b^{2s}a^{1+6t+3rs}, a^{3j}\}\cup \{a^{-1}, b^{2s}a^{-1-3t+3s}, b^sa^{-1-6t-3rs+6s}, a^{-3j}\}$, and the automorphism of $G$ induced by $a\mapsto a^{-1}$, $b\mapsto ba^{-3s^{-1}t-3r+6}$, fixes $S$ setwise, contrary to Proposition~\ref{prop=half-arc-transitive no involution}.

Let $i\neq 0$. Since the map $a\mapsto a$, $b\mapsto ba^{3m}$ induces an automorphism of $G$, we may further assume that $a,b^i\in R$ and $(b^i)^\sigma=b^i$, forcing that $b^\sigma=b$. In this case, $\sigma$ is induced by $a\mapsto b^sa^{1+3t}$, $b\mapsto b$. In particular, $R=\{a, b^sa^{1+3t}, b^{2s}a^{1+6t}, b^i\}$ and $S=R\cup R^{-1}=\{a, a^{-1}, b, b^{-1}, b^sa^{1+3t}, b^{2s}a^{-1-3t+3s}, b^{2s}a^{1+6t}, b^sa^{-1-6t+6s}\}$, where $s,t=0,1,2$. For $s=0$, we have $S=\{a, a^{-1}, b, b^{-1}, a^{1+3t}, a^{-1-3t}, a^{1+6t}, a^{-1-6t}\}$, and  the automorphism of $G$ induced by $a\mapsto a^{-1}$, $b\mapsto b$, fixes $S$ setwise, a contradiction. For $s=1$ or $2$, we have $S=\{a, a^{-1}, b, b^{-1}, ba^{1+3t}, b^2a^{-1-3(t-1)}, b^2a^{1+6t}, ba^{-1-6(t-1)}\}$, where $t=0,1,2$.
If $t=2$ then $S=\{a, a^{-1}, b, b^{-1}, ba^7, b^2a^5, b^2a^4, ba^2\}$, and  the automorphism of $G$ induced by $a\mapsto a^{-1}$, $b\mapsto b$, fixes $S$ setwise, a contradiction. If $t=0$ or $1$ then $S=\{a, a^{-1}, b, b^{-1}, ba, b^2a^2, b^2a, ba^5\}$ or $\{a, a^{-1}, b, b^{-1}, ba^4, b^2a^{-1}, b^2a^7, ba^{-1}\}$, but in both cases, by Magma~\cite{magma} $\Gamma$ is not edge-transitive, a contradiction.

\medskip
\f{\bf Case 2:} $G=G_2(3)=\langle a,b,c\ |\ a^3=b^3=c^3=1, [a,b]=c, [a,c]=[b,c]=1 \rangle$.

By Observation, $a,b\in R$. Let us consider the two cases depending whether the fixed vertex of $\s$ in $R$ lines in the center $Z(G)$ of $G$. Note that $Z(G)=\langle c\rangle$.

Assume the the fixed vertex of $\s$ in $R$ belongs to $Z(G)$. Then $c^\s=c$, and $a,b$ lines in one orbit of $\s$. If necessary, replace $\s$ by $\s^2$, and then we may assume that $\sigma$ is induced by $a\mapsto b$, $b\mapsto a^{-1}b^jc^k$, $c\mapsto c$. Since $\s^3=1$, we have $a=(a^{-1}b^jc^k)^\s=a^{-j}b^{j^2-1}c^{kj+k-j+2^{-1}j^2(j-1)}$, which implies that $j=-1$. Thus, $R=\{a, a^\sigma, a^{\sigma^2}, c^n\}=\{a,b,a^{-1}b^{-1}c^{k},c^n\}$ and $S=\{a,b,a^{-1}b^{-1}c^{k},c^n\}\cup \{a^{-1}, b^{-1}, abc^{-k-1}, c^{-n}\}$. The automorphism of $G$ induced by $a\mapsto a^{-1}$, $b\mapsto abc^{-k-1}$, $c\mapsto c^{-1}$, fixes $S$ setwise, contrary to Proposition~\ref{prop=half-arc-transitive no involution}.

Assume the the fixed vertex of $\s$ in $R$ is not in $Z(G)$. Noting that $\Aut(G)$ is transitive on elements  in $G\setminus Z(G)$, we may assume that $a^\s=a$. Then $\sigma$ is induced by $a\mapsto a$, $b\mapsto a^ib^jc^k$, $c\mapsto c^j$, where  $i,j,k\in \mathbb{Z}_3$. Since $\s^3=1$, we have  $b=b^{\sigma^3}=(a^ib^jc^k)^{\s^2}=a^sb^{j^3}c^t$ for some $s,t\in \mz_3$, and hence $j^3=1$, implying that $j=1$.
It follows that $R=\{a, b, b^\sigma, b^{\sigma^2}\}=\{a,b,a^ibc^k,a^{2i}bc^{2k}\}$ and $S=\{a,b,a^ibc^k,a^{2i}bc^{2k}\}\cup \{a^{-1}, b^{-1}, a^{-i}b^{-1}c^{-k-i}, a^{-2i}b^{-1}c^{-2k-2i}\}$. This is impossible because the automorphism of $G$ induced by $a\mapsto a^{-1}$, $b\mapsto b$, $c\mapsto c^{-1}$, fixes $S$ setwise.\hfill\qed

Now we are ready to prove Theorem~\ref{maintheorem}.

\medskip
\f{\bf Proof of Theorem~\ref{maintheorem}:} Let $\Gamma$ be a half-arc-transitive graph of order $p^3$ and valency $6$ or $8$, and let $A=\Aut(\Gamma)$.
By Observation, $\Gamma=\Cay(G,S)$ and $\overrightarrow{\Gamma}=\Cay(G,R)$ for some group $G=G_1(p)$ or $G_2(p)$ with $p\geq 3$, where $G=\langle R\rangle$ and $S=R\cup R^{-1}$. By Lemma~\ref{lem=no half 27 verties}, $p\geq 5$, and by the half-arc-transitivity of $\Gamma$, $A=\Aut(\overrightarrow{\Gamma})$ and $A_1$ is transitive on $R$.
Since $G=\langle R\rangle$, $\Aut(G,R)$ acts faithfully on $R$, and since $|R|<5$, $\Aut(G,R)$ is a $p'$-group. Since $G$ is a non-abelian group of order $p^3$, $G$ is $2$-genetic, that is, each normal subgroup of $G$ can be generated by two elements. By Proposition~\ref{3,5,7,11}, $\overrightarrow{\Gamma}$ is normal, and by Proposition~\ref{prop=Cayleynormal}, $A_1=\Aut(G,R)$. Since $A=\Aut(\Gamma)=\Aut(\overrightarrow{\Gamma})$, we have $A_1=\Aut(G,S)$, and so $\Gamma$ is normal.

The theorem is true for $G=G_1(p)$ by Proposition~\ref{prop=metacirculant graph}. Now assume $G=G_2(p)$. If $\Gamma$ has valency $8$, the theorem is also true by Lemma~\ref{lemma=necessity valency 8}. We may further assume that $\Gamma$ has valency $6$, that is, $|R|=3$. To finish the proof, it suffices to drive contradictions.

By Observation, $R=\{a,b,a^ib^jc^k\}$, where $i,j,k\in\mz_p$. Since $\Aut(G,R)$ is transitive on $R$, there exists $\alpha\in \Aut(G_2(p))$ of order $3$ permuting the elements in $R$ cyclically. If necessary, replace $\a$ by $\a^2$, and then we may assume that $\a$ is induced by
$a\mapsto b$, $b\mapsto a^ib^jc^k$, $c\mapsto c^{-i}$. Thus, $a=(a^ib^jc^k)^\alpha=a^{ij}b^{i+j^2}c^{-i^2j-2^{-1}ij^2(j-1)+k(j-i)}$, and so we have:

\parbox{8cm}{
\begin{eqnarray*}
&& ij=1;\\
&&i+j^2=0;\\
&&-i^2j+k(j-i)-2^{-1}ij^2(j-1)=0.
\end{eqnarray*}}\hfill
\parbox{1cm}{
\begin{eqnarray}
\label{eq1}\\ \label{eq2}\\\label{eq3}
\end{eqnarray}}

By Eqs~(\ref{eq1}) and (\ref{eq2}), $j^3+1=0$, implying  $(j+1)(j^2-j+1)=0$. Thus, either $j+1=0$ or $j^2-j+1=0$. If $j+1=0$ then $j=-1$. By Eq~(\ref{eq2}), $i=-1$ and so $S=\{a,b,a^{-1}b^{-1}c^k\}\cup \{a^{-1},b^{-1},abc^{-1-k}\}$, but the automorphism of $G$ induced by $a \mapsto a^{-1}$, $b \mapsto abc^{-1-k}$, $c \mapsto c^{-1}$, fixes $S$ setwise, contrary to Proposition~\ref{prop=half-arc-transitive no involution}.
If $j^2-j+1=0$ then by Eq~(\ref{eq2}), $i=1-j$, and since $ij=1$, Eq~(\ref{eq3}) implies that $-i+k(j-i)-2^{-1}j(j-1)=j-1+k(j+j-1)+2^{-1}=(2j-1)(k+2^{-1})=0$. It follows that either $2j-1=0$ or $k+2^{-1}=0$.
For $2j-1=0$, we have $j=2^{-1}$ and $i=1-j=1-2^{-1}=2^{-1}$, but Eq~(\ref{eq1}) gives rise to the contradiction that $p=3$.
For $k+2^{-1}=0$, we have  $S=\{a,b,a^{1-j}b^{j}c^{-2^{-1}}\}\cup \{a^{-1},b^{-1},a^{j-1}b^{-j}c^{-2^{-1}}\}$, and the automorphism of $G$ induced by $a \mapsto a^{-1}$, $b \mapsto b^{-1}$, $c \mapsto c$, fixes $S$ setwise, a contradiction. \hfill\qed

\medskip
\f {\bf Acknowledgements:} This work was supported by the National Natural Science Foundation of China (11231008, 11571035).


\begin{thebibliography}{99}

\bibitem{Alspach1}
B. Alspach, D. Maru\v{s}i\v{c} and L. Nowitz, Constructing graphs which are $1/2$-transitive, {\em J. Austral. Math. Soc. Ser. A} {\bf56} (1994), 391-402.
\bibitem{Alspach}
B. Alspach and M.Y. Xu, 1/2-transitive graphs of order $3p$, {\em J. Alg. Combin}. {\bf3} (1994), 347-355.
\bibitem{magma}
W. Bosma, C. Canon and C. Playoust, The MAGMA algebra system I: The user language, {\em J. Symbolic Comput}. {\bf24} (1997), 235-265.
\bibitem{Bouwer}
I.Z. Bouwer, Vertex and edge transitive but not 1-transitive graphs, {\em Canad. Math. Bull}. {\bf13} (1970), 231-237.
\bibitem{Chao}
C.Y. Chao, On the classification of symmetric graphs with a prime number of vertices, {\em Trans. Am. Math. Soc}. {\bf158} (1971), 247-256.
\bibitem{Cheng and Oxley}
Y. Cheng and J. Oxley, On weakly symmetric graphs of order twice a prime, {\em J. Combin. Theory Ser. B} {\bf42} (1987), 196-211.
\bibitem{Conder nonabelian stabilizer}
M.D.E. Conder and D. Maru\v{s}i\v{c}, A tetravalent half-arc-transitive graph with nonabelian vertex stabilizer, {\em J. Combin. Theory Ser. B} {\bf88} (2003), 67-76.
\bibitem{Conder recent discoverry}
M.D.E. Conder, P. Poto\v{c}nik and P. \v{S}parl, Some recent discoveries about half-arc-transitive graphs, {\em Ars Math. Contemp}. {\bf8} (2015), 149-162.
\bibitem{Dushaofei}
S.F. Du and M.Y. Xu, Vertex-primitive $1/2$-arc-transitive graphs of smallest order, {\em Comm. Algebra} {\bf27} (1999), 163-171.
\bibitem{Feng-valency 4}
Y.Q. Feng, J.H. Kwak and X.Y. Wang, Tetravalent half-arc-transitive graphs of order $2pq$, {\em J. Alg. Combin}. {\bf33} (2011), 543-553.
\bibitem{Feng}
Y.Q. Feng, J.H. Kwak, M.Y. Xu and J.X. Zhou, Tetravalent half-arc-transitive graphs of order $p^4$, {\em European J. Combin}. {\bf29} (2008), 555-567.


\bibitem{Godsil}
C.D. Godsil, On the full automorphism group of a graph, {\em Combinatorica} {\bf1} (1981), 243-256.
\bibitem{Holt}
D.F. Holt, A graph which is edge transitive but not arc transitive, {\em J. Graph Theory} {\bf5} (1981), 201-204.
\bibitem{Huppert}
B. Huppert, {\em Endliche Gruppen I}, Springer, Verlag, 1979.
\bibitem{Kutnar}
K. Kutnar, D. Maru\v{s}i\v{c}, P. \v{S}parl, R.J. Wang and M.Y. Xu, Classification of half-arc-transitive graphs of order $4p$, {\em European J. Combin}. {\bf34} (2013), 1158-1176.
\bibitem{Li1}
C.H. Li, On isomorphism of connected Cayley graphs, {\em Discrete Math}. {\bf178} (1998), 109-122.
\bibitem{Li2}
C.H. Li and H.S. Sim, On half-transitive metacirculant graphs of prime-power order, {\em J. Combin. Theory Ser. B} {\bf81} (2001), 45-57.

\bibitem{Marusic group action}
D. Maru\v{s}i\v{c}, Half-transitive group actions on finite graphs of valency $4$, {\em J. Combin. Theory Ser. B} {\bf 73} (1998), 41-76.
\bibitem{Marusic large stabilizer}
D. Maru\v{s}i\v{c}, Quartic half-arc-transitive graphs with large vertex stabilizers, {\em Discrete Math}. {\bf299} (2005), 180-193.
\bibitem{Marusic}
D. Maru\v{s}i\v{c}, Vertex transitive graphs and digraphs of order $p^k$, {\em Ann. Discrete Math}. {\bf115} (1985), 115-128.
\bibitem{Marusic and Praeger}
D. Maru\v{s}i\v{c} and C.E. Praeger, Tetravalent graphs admitting half-transitive group actions: alternating cycles, {\em J. Combin. Theory Ser. B} {\bf75} (1999), 188-205.
\bibitem{Dragan-valency4}
D. Maru\v{s}i\v{c} and P. Sparl, On quartic half-arc-transitive metacirculants, {\em J. Alg. Combin}. {\bf28} (2008), 365-395.

\bibitem{Primoz-valency4}
P. Sparl, On the classification of quartic half-arc-transitive metacirculants, {\em Discrete Math}. {\bf309} (2009), 2271-2283.
\bibitem{Robinson}
D.J. Robinson, {\em A course in the theory of groups}, second edition, Springer, New York, 1996.
\bibitem{Tutte}
W.T. Tutte, {\em Connectivity in Graphs}, University of Toronto Press, Toronto, 1966.
\bibitem{RJWang}
R.J. Wang, Half-transitive graphs of order a product of two distinct primes, {\em Comm. Algebra} {\bf22} (1994), 915-927.
\bibitem{Wang}
Y. Wang, Y.Q. Feng and J.X. Zhou, Cayley digraphs of $2$-genetic groups of prime-power order, submitted.
\bibitem{Xu}
M.Y. Xu, Automorphism groups and isomorphisms of Cayley digraphs, {\em Discrete Math}. {\bf182} (1998), 309-319.
\bibitem{Xu1}
M.Y. Xu, Half-transitive graphs of prime-cube order, {\em J. Alg. Combin}. {\bf1} (1992), 275-282.



\end{thebibliography}
\end{document}